\documentclass{birkmult}

\theoremstyle{plain}
\newtheorem{theorem}{Theorem}[section]

\newtheorem{proposition}[theorem]{Proposition}
\newtheorem{lemma}[theorem]{Lemma}

\theoremstyle{definition}

\theoremstyle{remark}

\begin{document}

\title[Unitary systems and wavelet sets]{Unitary systems and wavelet sets}

\author[Larson]{David R.~Larson}

\address{Department of Mathematics\\
Texas A{\&}M University\\
College Station, TX 77843-3368, U.S.A.}

\email{larson@math.tamu.edu}

\begin{abstract} A wavelet is a special case of a vector in a separable Hilbert space
that generates a basis under the action of a  system of unitary
operators defined in terms of translation and dilation operations.
We will describe an operator-interpolation approach to wavelet
theory using the local commutant of a unitary system. This is an
application of the theory of operator algebras to wavelet theory.
The concrete applications to wavelet theory include results obtained
using specially constructed families of wavelet sets. The main
section of this paper is section 5, in which we introduce the
interpolation map $\sigma$ induced by a pair of wavelet sets, and
give an exposition of its properties and its utility in constructing
new wavelets from old. The earlier sections build up to this,
establishing terminology and giving examples. The main theoretical
result is the Coefficient Criterion, which is described in Section
5.2.2, and which gives a matrix valued function criterion
specificing precisely when a function with frequency support
contained in the union of an interpolation family of wavelet sets is
in fact a wavelet. This can be used to derive Meyer's famous class
of wavelets using an interpolation pair of Shannon-type wavelet sets
as a starting point. Section 5.3 contains a new result on
interpolation pairs of wavelet sets: a proof that every pair of sets
in the generalized Journe family of wavelet sets is an interpolation
pair.  We will discuss some results that are due to this speaker and
his former and current students. And we finish in section 6 with a
discussion of some open problems on wavelets and frame-wavelets.
\end{abstract}


\thanks{This research was supported by a grant from the NSF.}

\maketitle


\section{Introduction}

  A wavelet is a special case of a
vector
 in a separable Hilbert space that generates a basis under the
 action of a collection, or "system", of unitary operators defined in terms of
 translation and dilation operations.  This approach to wavelet theory goes
 back, in particular, to earlier work of  Goodman, Lee and Tang [10] in the
 context of multiresolution analysis.  We will begin by describing the
 operator-interpolation approach to wavelet theory using the local commutant of
a
 system that was worked out by the speaker and his collaborators a few years
ago.  This
 is really an abstract application of the theory of operator algebras,
 mainly von Neumann algebras, to wavelet theory.  The concrete applications of
 operator-interpolation to wavelet theory include results obtained using
specially
 constructed families of wavelet sets.  In fact X. Dai and the speaker had
originally
 developed our theory of wavelet sets [5] specifically to take advantage of
their
 natural and elegant relationships with these wavelet unitary systems. We will
also discuss
 some new results and open questions.

The main idea in \textit{operator-theoretic interpolation} of
wavelets (and frames) is that  new wavelets can be obtained as
linear combinations of known ones using \textit{coefficients} which
are \textit{operators} (in fact, \emph{Fourier multipliers}) in a
certain class.  Both the ideas  and the essential computations
extend naturally to more general unitary systems and \emph{wandering
vectors}. Many of the methods work for more involved systems that
are important to applied harmonic analysis, such as Gabor and
generalized Gabor systems, and various types of \textit{frame}
unitary systems.

\subsection{Terminology}

 The set of all bounded linear operators on a Hilbert
space $H$ will be denoted by $B(H)$.  A \textit{bilateral shift} $U$
on $H$ is a unitary operator $U$ for which there exists a closed
linear subspace $E \subset H$ with the property that the family of
subspaces $\{U^nE: n \in \mathbb{Z}\}$ are orthogonal and give a
direct-sum decomposition of $H$.  The subspace $E$ is called a
\emph{complete wandering subspace} for $U$, and the
\textit{multiplicity} of $U$ is defined to be the dimension of $E$.
The \emph{strong operator topology} on $B(H$ is the topology of
pointwise convergence, and the \emph{weak operator topology} is the
weakest topology such that the vector functionals $\omega_{x,y}$ on
$B(H)$ defined by $A \mapsto \langle Ax, y \rangle$, $A \in B(H)$,
$x,y \in H$, are all continuous.  An \emph{algebra of operators} is
a linear subspace of $B(H)$ which is closed under multiplication. An
\emph{operator algebra} is an algebra of operators which is
\emph{norm-closed}. A subset $\mathcal{S} \subset B(H)$ is called
\emph{selfadjoint} if whenever $A \in \mathcal{S}$ then also $A^*
\in \mathcal{S}$.  A $C^*$-\emph{algebra} is a self-adjoint operator
algebra. A \emph{von Neumann algebra} is a $C^*$-algebra which is
closed in the weak operator topology.  For a unital operator
algebra, it is well known that being closed in the weak operator
topology is equivalent to being closed in the closed in the strong
operator topology. The \emph{commutant} of a set $\mathcal{S}$ of
operators in $B(H)$ is the family of all operators in $B(H)$ that
\emph{commute} with every operator in $\mathcal{S}$.  It is closed
under addition and multiplication, so is an algebra.  And it is
clearly closed in both the weak operator topology and the strong
operator topology.  We use the standard \emph{prime} notation for
the commutant. So the commutant of a subset $\mathcal{S} \subset
B(H)$ is denoted:~~ $\mathcal{S}^\prime := \{A \in B(H): AS = SA, ~~
S \in \mathcal{S}\}$. The commutant of a selfadjoint set of
operators is clearly a von Neumann algebra.  Moreover, by a famous
theorem of Fuglede every operator which commutes with a normal
operator $N$ also commutes with its adjoint $N^*$, and hence the
commutant of any set of \emph{normal} operators is also a von
Neumann algebra.  So, of particular relevance to this work, the
commutant of any set of \emph{unitary} operators is a von Neumann
algebra.

\subsubsection{Frames and Operators}
A sequence of vectors $\{f_j\}$ in a separable Hilbert space $ H$ is
a {\it frame} (or {\it frame sequence}) if there exist constants
$C_1, C_2
>0$ such that
$$
   C_1 \|f\|_2^2 \leq \sum_j |<f, f_j>|^2 \leq
C_2 \|f\|_2^2
$$
for all $f\in H$. If $C_1=C_2$ the frame is called \emph{tight}, and
if  $C_1=C_2 = 1$, $\{f_j\}$ is called a {\emph Parseval} frame.
(The term \emph{normalized tight} has also been used for this (cf
[14]). A vector $\xi$ is called a {\emph frame vector} for a unitary
system ${\mathcal U}$ if the set of vectors ${\mathcal U}\xi$ is a
frame for $H$.

A Riesz basis for a Hilbert space is a bounded unconditional basis.
Frames sequences are generalizations of Riesz bases. A number of the
basic aspects of a geometric, or operator-theoretic, approach to
discrete frame theory on Hilbert space arises from the fact that a
frame sequence is simply an "inner" direct summand of a Riesz basis.
The basic principle is that a Hilbert space frame sequence can be
dilated to a Riesz basis for larger Hilbert space.  We call this the
\emph{Frame Dilation Theorem}. In other words, for a given frame
sequence there is a larger Hilbert space and a Riesz basis for the
larger space such that the orthogonal projection from the larger
space onto the smaller space compresses the Riesz basis to the frame
sequence. We proved this at the beginning of [14], and used it to
prove the other results [14], and subsequently to prove some
applications to Hilbert C*-module theory jointly with M. Frank.  We
proved it first for Parseval frames, and then for general frames.
(We remark that this type of dilation result for frames was also
independently known and used independently by several others in
different contexts.)

It is interesting to note that the Parseval frame case of the Frame
Dilation Theorem can be derived easily from the purely atomic case
of a well known theorem of Naimark on projection valued measures. We
thank Chandler Davis for pointing this out to us at the Canadian
Operator Algebras Symposium in 1999. We (Han and I) basically proved
this special case of Naimark's theorem implicitely in the first
section of [14] without recognizing it was a special case of
Naimark's theorem, and then we proved the appropriate generalization
we needed for general (non-tight) frames.  Naimark's Dilation
Theorem basically states that a suitable positive operator valued
measure on a Hilbert space \emph{dilates} to a projection valued
measure on a larger Hilbert space.  That is, there is a projection
valued measure [PVM] on a larger Hilbert space such that the
orthogonal projection from the larger space unto the smaller space
compresses the [PVM] to the [POVM]. In the discrete (i.e. purely
atomic measure) case, it can be interpreted as stating that a
suitable sequence of positive operators dilates to a sequence of
projections.  The dilation theorem for a Parseval frame follows
easily from Naimark's Theorem applied to the [POVM] obtained by
replacing each vector $x_i$ in the frame sequence with the
elementary tensor operator $x_i \otimes x_i$, obtaining the atoms
for a [POVM] defined on all subsets of the index set for the frame.
The dilation theorem for a general (non-tight) frame does not seem
to follow directly from Naimark's theorem -- but it may follow from
a generalization of it.  We remark that some other generalizations
of the frame dilation theorem have been recently worked out, notably
by W. Czaja.

\subsubsection{Unitary Systems, Wandering Vectors, and Frame Vectors}

We define a \textit{unitary system} to be simply a countable
collection of unitary operators $\mathcal U$ acting on a Hilbert
space $H$ which contains the identity operator.  The
\textit{interesting} unitary systems all have additional structural
properties of various types.  (For instance, wavelet systems and
Gabor systems are both "ordered products" of two abelian groups: the
dilation and translation groups in the wavelet case, and the
modulation and translation groups in the Gabor case.)
  We will say that a vector $\psi \in H$ is \textit{wandering} for
$\mathcal U$ if the set
\begin{equation} \mathcal U \psi := \{U\psi : U \in \mathcal U\} \end{equation}
is an orthonormal set, and we will call $\psi$ a \textit{complete
wandering vector} for $\mathcal U$ if $\mathcal U \psi$ spans $H$.
This (abstract) point of view can be useful.  Write $\mathcal
W(\mathcal U)$ for the set of complete wandering vectors for
$\mathcal U$.  Further, a \emph{Riesz vector}  for $\mathcal U$ is a
vector $\psi$ such that $\mathcal U \psi$ is a Riesz basis for $H$
(indexed by the elements of $\mathcal U$), and a \emph{frame vector}
is a vector $\psi$ such that $\mathcal U \psi$ is a frame sequence
for $H$ (again using $\mathcal U$ as its index set), and we adopt
similar terminology for Parseval frame vectors and Bessel vectors.
We use $\mathcal RW(\mathcal U)$, $\mathcal F(\mathcal U)$,
$\mathcal PF(\mathcal U)$, $\mathcal B(\mathcal U)$ to denote,
respectively, the sets of Riesz vectors, frame vectors, Parseval
frame vectors, and Bessel vectors for $\mathcal U$.

One of the main tools in this work is the \emph{local commutant} of
a system of unitary operators (see section 3.2).  This is a natural
generalization of the commutant of the system, and like the
commutant it is a linear space of operators which is closed in the
weak and the strong operator topologies, but unlike the commutant it
is usually not selfadjoint, and is usually not closed under
multiplication.  It contains the commutant of the system, but can be
much larger than the commutant.  The local commutant of a wavelet
unitary system captures all the information about the wavelet system
in an essential way, and this gives the \emph{flavor} of our
approach to the subject.

\subsubsection{Normalizers}

If $U$ is a unitary operator and $\mathcal{A}$ is an operator
algebra, then $U$ is said to \emph{normalize} $\mathcal{A}$ if
~~$U^\star \cdot \mathcal{A} \cdot U = \mathcal{A}$~.  In the most
interesting cases of operator-theoretic interpolation: that is, for
those cases that yield the strongest structural results, the
relevant unitaries in the local commutant of the system normalize
the commutant of the system.

\subsection{Acknowledgement}

This article was written in response to an invitation by the
organizers of the 4th International Conference on Wavelet Analysis
and its Applications, Macau, China, December 2005 [WAA2005] to be a
keynote speaker.  We thank the organizers for their kind invitation
to present these notes.

\section{Wavelets}

For simplicity of presentation, much of the work in this article
will deal with one-dimensional wavelets, and in particular, the
dyadic case.  The other cases: non-dyadic wavelets and wavelets in
higher dimensions,
 are at least notationally
more complicated.

\subsection{One Dimension}

 A \textit{dyadic orthonormal} wavelet in one dimension is a unit
vector $\psi \in L^2(\mathbb{R}, \mu)$, with $\mu$ Lebesgue measure,
with the property that the set
\begin{equation} \{2^{\frac{n}{2}}\psi(2^n t - l) : n,l \in
\mathbb{Z}\}\end{equation} of all integral translates of $\psi$
followed by dilations by arbitrary integral powers of $2$, is an
orthonormal basis for $L^2(\mathbb{R},\mu)$.  The term \emph{dyadic}
refers to the dilation factor "$2$".  The term \emph{mother wavelet}
is also used in the literature for $\psi$. Then the functions
$$\psi_{n,l} := 2^{\frac{n}{2}} \psi(2^n t - l)$$ are called elements
of the wavelet basis generated by the "mother".  The functions
$\psi_{n,l}$ will not themselves be mother wavelets unless $n = 0$.

Let $T$ and $D$ be the translation (by $1$) and dilation (by $2$)
unitary operators in $B(L^2(\mathbb{R})$ given by $(Tf)(t) = f(t-1)$
and $(Df)(t) = \sqrt{2}f(2t)$.  Then
$$2^{\frac{n}{2}}\psi(2^nt - l) = (D^nT^l \psi)(t)$$
for all $n,l \in \mathbb{Z}$.  Operator-theoretically, the operators
$T, D$ are \textit{bilateral shifts} of \textit{infinite
multiplicity}.  It is obvious that $L^2([0,1])$, considered as a
subspace of $L^2(\mathbb{R})$, is a complete wandering subspace for
$T$, and that  $L^2([-2, -1] \cup [1, 2])$ is a complete wandering
subspace for $D$.

Let $\mathcal{U}_{D,T}$ be the unitary system defined by
\begin{equation} \mathcal{U}_{D,T}= \{D^nT^l : n,l \in \mathbb{Z}\}
\end{equation}
where $D$ and $T$ are the operators defined above.  Then $\psi$ is a
dyadic orthonormal wavelet if and only if $\psi$ is a complete
wandering vector for the unitary system $\mathcal{U}_{D,T}$. This
was our original motivation for developing the abstract unitary
system theory.  Write
\begin{equation} \mathcal{W}(D,T) := \mathcal{W}(\mathcal{U}_{D,T})
\end{equation}
  to denote the set of all dyadic orthonormal wavelets in one
  dimension.

An abstract interpretation is that, since $D$ is a bilateral shift
it has (many) complete wandering subspaces,  and a wavelet for the
system is a vector $\psi$ whose translation space (that is, the
closed linear span of $\{T^k: k \in \mathbb{Z}\}$ is a complete
wandering subspace for $D$.  Hence $\psi$ must generate an
orthonormal basis for the entire Hilbert space under the action of
the unitary system.

  In one dimension, there are non-dyadic orthonormal wavelets: i.e.
wavelets for all possible dilation factors besides $2$ (the dyadic
case). We said "possible", because the scales $\{0,1,-1\}$ are
excluded as scales because the dilation operators they would
introduce are not bilateral shifts. All other real numbers for
scales yield wavelet theories. In [5, Example 4.5 (x)] a family of
examples is given of three-interval wavelet sets (and hence
wavelets) for all scales $d \geq 2$, and it was noted there that
such a family also exists for dilation factors $1 < d \leq 2$. There
is some recent (yet unpublished) work that has been done, by REU
students and mentors, building on this, classifying finite-interval
wavelet sets for all possible real (positive and negative scale
factors).  I will mention this work, in passing, in my talk.

\subsection{N-Dimensions}

\subsubsection{Expansive Dilations}

Let $1 \leq m < \infty$, and let $A$ be an $n \times n$ real matrix
which is \textit{expansive} (equivalently, all (complex) eigenvalues
have modulus $>1$).  By a \emph{dilation - $A$ regular-translation
orthonormal wavelet} we mean a function $\psi \in L^2(\mathbb{R}^n)$
such that
\begin{equation} \{|det(A)|^{\frac{n}{2}} \psi(A^n t - (l_1, l_2, ..., l_n)^t) :
n,l \in \mathbb{Z}\} \end{equation}
where $t = (t_1, ..., t_n)^t$, is an orthonormal basis for
$L^2(\mathbb{R}^n ; m)$.  (Here $m$ is product Lebesgue measure, and
the superscript "t" means transpose.)

If $A \in M_n(\mathbb{R})$ is invertible (so in particular if $A$ is
expansive), then it is very easy to verify that the operator defined
by
\begin{equation} (D_Af)(t) = |det A|^{\frac12} f(At) \end{equation}
for $f \in L^2(\mathbb{R}^n)$, $t \in \mathbb{R}^n$, is
\emph{unitary}. For $1 \leq i \leq n$, let $T_i$ be the unitary
operator determined by translation by $1$ in the $i^{th}$ coordinate
direction.  The set (5) above is then
\begin{equation} \{D^k_A T^{l_1}_1 \cdot\cdot\cdot T^{l_n}_n \psi : k,l_i \in
\mathbb{Z}\} \end{equation}

If the dilation matrix $A$ is expansive, but the translations are
along some oblique lattice, then there is an invertible real $n
\times n$ matrix $T$ such that conjugation with $D_T$ takes the
entire wavelet system to a regular-translation expansive-dilation
matrix. This is easily worked out, and was shown in detail in [18]
in the context of working out a complete theory of unitary
equivalence of wavelet systems. Hence the wavelet theories are
equivalent.

\subsubsection{Non-Expansive Dilations}

Much work has been accomplished concerning the existence of wavelets
for dilation matrices $A$ which are not expansive.  Some of the
original work was accomplished in the Ph.D. theses of Q. Gu and  D.
Speegle, when they were together finishing up at Texas A\&M. Some
significant additional work was accomplished by Speegle and also by
others. In  [18], with Ionascu and Pearcy we proved that if an $n x
n$ real invertible matrix $A$ is not similar (in the nxn complex
matrices) to a unitary matrix, then the corresponding dilation
operator $D_A$ is in fact a bilateral shift of infinite
multiplicity.  If a dilation matrix were to admit any type of
wavelet (or frame-wavelet) theory, then it is well-known that a
necessary condition would be that the corresponding dilation
operator would have to be a bilateral shift of infinite
multiplicity.  I am happy to report that in very recent work [23],
with E. Schulz, D. Speegle, and K. Taylor, we have succeeded in
showing that this minimal condition is in fact sufficient: such a
matrix, with regular translation lattice, admits a (perhaps
infinite) tuple of functions, which collectively generates a
frame-wavelet under the action of this unitary system.

\section{More General Unitary Systems}

\subsection{Some Restrictions}

We note that \textit{most} unitary systems $\mathcal{U}$ do not have
complete wandering vectors.  For $\mathcal{W}(\mathcal{U})$ to be
nonempty, the set $\mathcal{U}$ must be very special. It must be
\emph{countable} if it acts separably (i.e. on a separable Hilbert
space), and it must be \emph{discrete} in the strong operator
topology because if $U,V \in \mathcal{U}$ and if $x$ is a wandering
vector for $\mathcal{U}$ then
$$\|U - V \| \geq \| Ux - Vx \| = \sqrt{2}$$
Certain other properties are forced on $\mathcal{U}$ by the presence
of a wandering vector.  (Or indeed, by the nontriviality of any of
the sets $\mathcal{W}(\mathcal{U})$, $\mathcal{RW}(\mathcal{U})$,
$\mathcal{F}(\mathcal{U})$, $\mathcal{PF}(\mathcal{U})$,
$\mathcal{B}(\mathcal{U})$.)  One purpose of [5] was to investigate
such properties. Indeed, it was a matter of some surprise to us to
discover that such a theory is viable even in some considerable
generality.  For perspective, it is useful to note that while
$\mathcal{U}_{D,T}$ has complete wandering vectors, the reversed
system $$\mathcal{U}_{T,D} = \{T^lD^n : n,l \in \mathbb{Z}\}$$
\textit{fails} to have a complete wandering vector. (A proof of this
was given in the introduction to [5].)

\subsection{The Local Commutant}

\subsubsection{ A Special Case: The System $\mathcal{U}_{D,T}$}

Computational aspects of operator theory can be introduced into the
wavelet framework in an elementary way.  Here is the way we
originally did it:  Fix a wavelet $\psi$ and consider the set of all
operators  $A \in  B(L^2(\mathbb{R}))$ which \textit{commute} with
the \emph{action} of dilation and translation on $\psi$.  That is,
require
\begin{equation}(A\psi)(2^nt-l) = A(\psi(2^nt - l))\end{equation}
 or equivalently
\begin{equation}D^nT^lA\psi = AD^nT^l\psi \end{equation} for all $n,l \in
\mathbb{Z}$.  Call this the \emph{local commutant of the wavelet
system $\mathcal{U}_{D,T}$ at the vector $\psi$}. (In our first
preliminary writings and talks we called it the \emph{point
commutant} of the system.) Formally, the local commutant of the
dyadic wavelet system on $L^2(\mathbb{R})$ is:
\begin{equation}\mathcal{C}_\psi(\mathcal{U}_{D,T}) := \{A \in
B(L^2(\mathbb{R})): (AD^nT^l - D^nT^lA)\psi = 0, \forall n,l \in
\mathbb{Z}\}\end{equation}
This is a linear subspace of $B(H)$ which is closed in the strong
operator topology, and in the weak operator topology, and it clearly
contains the \emph{commutant} of $\{ D,T \}$.

A motivating example is that if $\eta$ is any other wavelet, let $V
:= V_\psi^\eta$ be the unitary (we call it the \textit{interpolation
unitary}) that takes the basis $\psi_{n,l}$ to the basis
$\eta_{n,l}$.  That is, $V\psi_{n,l} = \eta_{n,l}$ for all $n,l \in
\mathbb{Z}$.  Then $\eta = V\psi$, so $VD^nT^l\psi = D^nT^lV\psi$
hence $V \in \mathcal{C}_\psi(\mathcal{U}_{D,T})$.

In the case of a pair of complete wandering vectors $\psi,\eta$ for
a general unitary system $\mathcal{U}$, we will use the same
notation $V_\psi^\eta$ for the unitary that takes the vector $U\psi$
to $U\eta$ for all $U \in \mathcal{U}$.

This simple-minded idea is reversible, so for every unitary $V$ in
$\mathcal{C}_\psi(\mathcal{U}_{D,T})$ the vector $V \psi$ is a
wavelet. This correspondence between unitaries in
$\mathcal{C}_\psi(D,T)$ and dyadic orthonormal wavelets is
one-to-one and onto (see Proposition ~3.1).   This turns out to be
useful, because it leads to some new formulas relating to
decomposition and factorization results for wavelets, making use of
the \textit{linear} and \textit{multiplicative} properties of
$\mathcal{C}_\psi(D,T)$.

It turns out (a proof is required) that the entire local commutant
of the system $\mathcal{U}_{D,T}$ at a wavelet $\psi$ is \emph{not}
closed under multiplication, but it also turns out (also via a
proof) that for \emph{most} (and perhaps \emph{all}) wavelets $\psi$
the local commutant at $\psi$ contains many noncommutative operator
algebras (in fact von Neumann algebras) as subsets, and their
unitary groups \emph{parameterize} norm-arcwise-connected families
of wavelets. Moreover, $\mathcal{C}_\psi(D,T)$ is closed under
\emph{left multiplication} by the commutant $\{ D,T \}^\prime$,
which turns out to be an abelian nonatomic von Neumann algebra.  The
fact that $\mathcal{C}_\psi(D,T)$ is a \emph{left module} under $\{
D,T \}^\prime$ leads to a method of obtaining new wavelets from old,
and of obtaining connectedness results for wavelets, which we called
\textit{operator-theoretic interpolation} of wavelets in [5], (or
simply \emph{operator-interpolation}).

\subsubsection{General Systems}

More generally, let $\mathcal{S} \subset B(H)$ be a set of
operators, where $H$ is a separable Hilbert space, and let $x \in H$
be a nonzero vector, and \textit{formally} define the \textit{local
commutant} of $\mathcal{S}$ at $x$ by
$$\mathcal{C}_x(\mathcal{S}) := \{A \in B(H) : (AS - SA)x = 0, S \in
\mathcal{S}\}$$

As in the wavelet case, this is a weakly and strongly closed linear
subspace of $B(H)$ which contains the commutant $\mathcal{S}^\prime$
of $\mathcal{S}$.  If $x$ is \textit{cyclic} for $\mathcal{S}$ in
the sense that span$(\mathcal{S}x)$ is dense in $H$, then $x$
\textit{separates} $\mathcal{C}_x(\mathcal{S})$ in the sense that
for $S \in \mathcal{C}_x(\mathcal{S})$, we have $Sx = 0$ iff $x =
0$.
Indeed, if $A \in \mathcal{C}_x(\mathcal{S})$ and if $Ax = 0$, then
for any $S \in \mathcal{S}$ we have $ASx = SAx = 0$, so $A
\mathcal{S}x = 0$, and hence $A = 0$.

If $A \in \mathcal{C}_x(\mathcal{S})$ and $B \in
\mathcal{S}^{\prime}$,  let $C = BA$.  Then for all $S \in
\mathcal{S}$,
$$(CS - SC)x = B(AS)x - (SB)Ax = B(SA)x - (BS)Ax = 0$$
because $ASx = SAx$ since $A \in \mathcal{C}_x(\mathcal{S})$, and
$SB = BS$ since $B \in \mathcal{S}^{\prime}$.
Hence $\mathcal{C}_x(\mathcal{S})$ is closed under left
multiplication by operators in $\mathcal{S}^{\prime}$. That is,
$\mathcal{C}_x(\mathcal{S}$ is a \emph{left module} over
$\mathcal{S}^{\prime}$.

It is interesting that, if in addition $\mathcal{S}$ is a
multiplicative semigroup, then in fact $\mathcal{C}_x(\mathcal{S})$
is identical with the commutant
$\mathcal{S}^\prime$ so in this case the commutant is not a new
structure.  To see this, suppose $A \in \mathcal{C}_x(\mathcal{S})$.
Then for each $S,T \in \mathcal{S}$ we have $ST \in \mathcal{S}$,
and so
$$AS(Tx) = (ST)Ax = S(ATx) = (S)Tx$$
So since $T \in \mathcal{S}$ was arbitrary and span$(\mathcal{S}x) =
H$, it follows that $AS = SA$.

\begin{proposition} If $\mathcal{U}$ is any unitary system for which
$\mathcal{W}(\mathcal{U}) \neq \emptyset$, then for any $\psi \in
\mathcal{W}(\mathcal{U})$
$$\mathcal{W}(\mathcal{U}) = \{U\psi : U
\textit{ is a unitary operator in }  \mathcal{C}_\psi(\mathcal{U})\}
$$ and the correspondence $U \rightarrow U\psi$ is one-to-one.
\end{proposition}

\begin{proposition} Let $\mathcal{U}$ be a unitary system on a Hilbert space
$H$.  If $\psi$ is a complete wandering vector for $\mathcal{U}$ ,
then:
\begin{itemize}
\item[(i)] $\mathcal{RW(U)} = \{A\psi : A \textit{ is an operator in }
\mathcal{C}_\psi(\mathcal{U}) \textit{ that is invertible in }
B(H)\};$
\item[(ii)] $\mathcal{F(U)} = \{A\psi : A \textit{ is an operator in }
\mathcal{C}_\psi(\mathcal{U}) \textit{ that is surjective}\};$
\item[(ii)] $\mathcal{PF(U)} = \{A\psi : A \textit{ is an operator in }
\mathcal{C}_\psi(\mathcal{U}) \textit{ that is a co-isometry}\};$
\item[(ii)] $\mathcal{B(U)} = \{A\psi : A \textit{ is an operator in }
\mathcal{C}_\psi(\mathcal{U}) \}$
\end{itemize}

\end{proposition}

\subsection{Operator-Theoretic Interpolation }

Now suppose $\mathcal{U}$ is a unitary system, such as
$\mathcal{U}_{D,T}$, and suppose $\{\psi_1, \psi_2, \dots, \psi_m \}
\subset \mathcal{W}(\mathcal{U})$.  (In the case of
$\mathcal{U}_{D,T}$, this means that $(\psi_1, \psi_2, \dots, \psi_n
)$ is an n-tuple of wavelets.

Let $(A_1, A_2, \dots, A_n)$ be an n-tuple of operators in the
commutant  $\mathcal{U}^{\prime}$ of $\mathcal{U}$, and let $\eta$
be the vector
$$ \eta := A_1 \psi_1 + A_2 \psi_2 + \dots + A_n \psi_n ~.$$
Then

$$ \eta = A_1 \psi_1 + A_2 V_{\psi_1}^{\psi_2}\psi_1 + \dots A_n
V_{\psi_1}^{\psi_n}\psi_1 $$
\begin{equation} ~~= (A_1 + A_2
V_{\psi_1}^{\psi_2} + \dots + A_n V_{\psi_1}^{\psi_n}) \psi_1~~.
\end{equation}

We say that $\eta$ is obtained by $\emph{operator interpolation}$
from $\{\psi_1, \psi_2, \dots, \psi_m \}$.  Since
$\mathcal{C}_{\psi_1}(\mathcal{U})$ is a left $\mathcal{U}^{\prime}$
- module, it follows that the operator
\begin{equation} A := A_1 + A_2 V_{\psi_1}^{\psi_2} + \dots A_n
V_{\psi_1}^{\psi_n}
\end{equation}
is an element of $\mathcal{C}_{\psi_1}(\mathcal{U})$.  Moreover, if
$B$ is another element of $\mathcal{C}_{\psi_1}(\mathcal{U})$ such
that $\eta = B \psi_1$, then $~A - B ~\in
\mathcal{C}_{\psi_1}(\mathcal{U})$ and $(A - B)\psi_1 = 0$.  So
since $\psi_1$ separates $\mathcal{C}_{\psi_1}(\mathcal{U})$ it
follows that $A = B$.  Thus $A$ is the \emph{unique} element of
$\mathcal{C}_{\psi_1}(\mathcal{U})$ that takes $\psi_1$ to $\eta$.
Let $\mathcal{S}_{\psi_1, \dots, \psi_n}$ be the family of all
finite sums of the form
$$ \sum^n_{i=0} A_i V_{\psi_1}^{\psi_i} ~~~.$$
This is the left module of $\mathcal{U}^{\prime}$ generated by $\{I,
V_{\psi_1}^{\psi_2}, \dots, V_{\psi_1}^{\psi_n}\}$.  It is the
$\mathcal{U}^{\prime}$-\emph{linear span} of $\{I,
V_{\psi_1}^{\psi_2}, \dots, V_{\psi_1}^{\psi_n}\}$.

Let
\begin{equation}
\mathcal{M}_{\psi_1, \dots, \psi_n} := (\mathcal{S}_{\psi_1, \dots,
\psi_n}) \psi_1
\end{equation}

\noindent So
\[ \mathcal{M}_{\psi_1, \dots, \psi_n} ~=~ \left\{\sum^n_{i=0} A_i
\psi_i ~: ~A_i \in \mathcal{U}^{\prime}\right\} ~~~.\]

\noindent We call this the \emph{interpolation space} for
$\mathcal{U}$ generated by $(\psi_1, \dots, \psi_n)$.  From the
above discussion, it follows that for every vector $\eta \in
\mathcal{M}_{\psi_1, \psi_2, \dots, \psi_n}$ there exists a unique
operator $A \in \mathcal{C}_{\psi_1}(\mathcal{U})$ such that $\eta =
A \psi_1$, and moreover this $A$ is an element of
$\mathcal{S}_{\psi_1, \dots, \psi_n}$.

\subsubsection{Normalizing the Commutant}
In certain essential cases (and we are not sure how general this
type of case is) one can prove that an interpolation unitary
$V_\psi^\eta$ \emph{normalizes} the commutant $\mathcal{U}^{\prime}$
of the system in the sense that $V_\eta^\psi \mathcal{U}^{\prime}
V_\psi^\eta = \mathcal{U}^{\prime}$.  (Here, it is easily seen that
$(V_\psi^\eta)^* = V_\eta^\psi$.)  Write $V := V_\psi^\eta$.  If $V$
normalizes $\mathcal{U}^\prime$, then the algebra, before norm
closure, generated by
$\mathcal{U}^\prime$ and $V$ is the set of all finite sums (trig
polynomials) of the form $\sum A_n V^n$, with coefficients $A_n \in
\ \mathcal{U}^\prime$, $n \in \mathbb{Z}$. The closure in the strong
operator topology is a von Neumann algebra. Now suppose further that
{\em every power\/} of $V$ is contained in
$\mathcal{C}_\psi(\mathcal{U})$. This occurs only in special cases,
yet it occurs frequently enough to yield some general methods. Then
since $\mathcal{C}_\psi(\mathcal{U})$ is a SOT-closed linear
subspace which is closed under left multiplication by
$\mathcal{U}^\prime$, this von Neumann algebra is contained in
$\mathcal{C}_\psi(\mathcal{U})$, so its unitary group parameterizes
a norm-path-connected subset of $\mathcal{W}(\mathcal{U})$ that
contains $\psi$ and $\eta$ via the correspondence $U\to U\psi$.

In the special case of \emph{wavelets}, this is the basis for the
work that Dai and I did in [5, Chapter 5] on operator-theoretic
interpolation of wavelets.  In fact, we specialized there and
\emph{reserved} the term \emph{operator-theoretic interpolation} to
refer explicitly to the case when the interpolation unitaries
normalize the commutant. In some subsequent work, we \emph{loosened}
this restriction yielding our more general definition given in this
article, because there are cases of interest in which we weren't
able to prove normalization.  However, it turns out that if $\psi$
and $\eta$ are $s$-elementary wavelets (see section 4.4), then
indeed $V^\eta_\psi$ normalizes $\{D,T\}'$. (See Proposition 5.3.)
Moreover, $V^\eta_\psi$ has a very special form:\ after conjugating
with the Fourier transform, it is a composition operator with a
symbol $\sigma$ that is a natural and very computable
measure-preserving transformation of $\mathbb{R}$. In fact, it is
precisely this special form for $V^\eta_\psi$ that allows us to make
the computation that it normalizes $\{D,T\}'$. On the other hand, we
know of no pair $(\psi,\eta)$ of wavelets for which $V^\eta_\psi$
fails to normalize $\{D,T\}'$. The difficulty is simply that in
general it is very hard to do the computations.  This is stated as
Problem 2 in the final section on \emph{Open Problems}.

In the wavelet case $\mathcal{U}_{D,T}$ , if
$\psi \in \mathcal{W}(D,T)$ then it turns out that
$\mathcal{C}_{\psi}(\mathcal {U}_{D,T})$ is in fact $\emph{much
larger}$ than $(\mathcal{U}_{D,T})^{\prime} = \{D,T\}^{\prime}$ ,
underscoring the fact that $\mathcal{U}_{D,T}$ is NOT a group.  In
particular, $\{D,T \}^{\prime}$ is abelian while $\mathcal{C}_ \psi
(\mathcal{D,T})$ is nonabelian for every wavelet $\psi$.  (The proof
of these facts are contained in [5].)

\subsubsection{Interpolation Pairs of Wandering Vectors}

In some cases where a pair $\psi, \eta$ of vectors in
$\mathcal{W(U)}$ are given  it turns out that the unitary $V$ in
$\mathcal{C}_\psi(\mathcal{U})$ with $V\psi = \eta$ happens to be a
\emph{symmetry} (i.e. $V^2 = I$).  Such pairs are called
\emph{interpolation pairs} of wandering vectors, and in the case
where $\mathcal{U}$ is a wavelet system, they are called
interpolation pairs of wavelets. Interpolation pairs are more
prevalent in the theory, and in particular the wavelet theory, than
one might expect. In this case (and in more complex generalizations
of this) certain linear combinations of complete wandering vectors
are themselves complete wandering vectors -- not simply complete
Riesz vectors.

\begin{proposition} Let $\mathcal{U}$ be a unitary system, let
$\psi, \eta \in \mathcal{W(U)}$, and let $V$ be the unique operator
in $\mathcal{C}_\psi(\mathcal{U})$ with $V \psi = \eta$.  Suppose
$$V^2 = I.$$ Then
$$ \cos \alpha \cdot \psi \textit{ +  } i \sin \alpha \cdot \eta \in
\mathcal{W(U)}$$ for all $0 \leq \alpha \leq 2\pi$.

\end{proposition}

The above result can be thought of as the \emph{prototype} of our
operator-theoretic interpolation results.  It is the second most
elementary case.  (The most elementary case is described in the
context of the exposition of Problem 4 in the final section.) More
generally, the scalar $\alpha$ in Proposition ~3.3 can be replaced
with an appropriate \emph{self-adjoint operator} in the commutant of
$\mathcal{U}$.  In the wavelet case, after conjugating with the
Fourier transform, which is a unitary operator, this means that
$\alpha$ can be replaced with a wide class of nonnegative
dilation-periodic (see definition below) bounded measurable
functions on $\mathbb{R}$.

\section{Wavelet Sets}

Wavelet sets belong to the theory of wavelets via the Fourier
Transform.  We will do most of this section in a tutorial-style, to
make the concepts more accessible to students and colleagues who are
not already familiar with them.

\subsection{Fourier Transform}

We will use the following form of the Fourier--Plancherel transform
$\mathcal{F}$ on $\mathcal{H} = L^2(\mathbb{R})$, which is a form
that is \emph{normalized} so it is a unitary transformation, a
property that is desirable for our treatment.

If $f,g\in L^1(\mathbb{R}) \cap L^2(\mathbb{R})$ then
\begin{equation}\label{eq18}
(\mathcal{F}f)(s) := \frac1{\sqrt{2\pi}} \int_{\mathbb{R}} e^{-ist} f(t)dt :=
\hat f(s),
\end{equation}
and
\begin{equation}\label{eq19}
(\mathcal{F}^{-1}g)(t) = \frac1{\sqrt{2\pi}} \int_{\mathbb{R}} e^{ist}g(s)ds.
\end{equation}
We have
\[
(\mathcal{F} T_\alpha f)(s) = \frac1{\sqrt{2\pi}} \int_{\mathbb{R}} e^{-ist}
f(t-\alpha)dt = e^{-is\alpha} (\mathcal{F}f)(s).
\]
So $\mathcal{F} T_\alpha \mathcal{F}^{-1} g=e^{-is\alpha}g$. For
$A\in\mathcal{B}(\mathcal{H})$ let $\hat A$ denote
$\mathcal{F}A\mathcal{F}^{-1}$. Thus
\begin{equation}\label{eq20}
\widehat T_\alpha = M_{e^{-i\alpha s}},
\end{equation}
where for $h\in L^\infty$ we use $M_h$ to denote the multiplication operator
$f\to hf$. Since $\{M_{e^{-i\alpha s}}\colon \ \alpha\in\mathbb{R}\}$ generates
the m.a.s.a.\ $\mathcal{D}(\mathbb{R}) := \{M_h\colon \ h\in
L^\infty(\mathbb{R})\}$ as a von Neumann algebra, we have
\[
\mathcal{F}\mathcal{A}_T \mathcal{F}^{-1} = \mathcal{D}(\mathbb{R}).
\]
Similarly,
\begin{align*}
(\mathcal{F}D^nf)(s) &= \frac1{\sqrt{2\pi}} \int_{\mathbb{R}} e^{-ist} (\sqrt
2)^n f(2^nt)dt\\
&= (\sqrt 2)^{-n}\cdot \frac1{\sqrt{2\pi}} \int_{\mathbb{R}} e^{-i2^{-n}st}
f(t)dt\\
&= (\sqrt 2)^{-2} (\mathcal{F}f)(2^{2^{-n}s}) = (D^{-n}\mathcal{F}f)(s).
\end{align*}
So $\widehat D^n = D^{-n} = D^{*n}$. Therefore,
\begin{equation}\label{eq21}
\widehat D = D^{-1} = D^*.
\end{equation}

If $f$ is an $L^2(\mathbb{R})$ function, as usual we write $\widehat
f(s) = (\mathcal{F}(f))(s)$.  If $\rho(s)$ is a real-valued function
such that $\widehat f(s) = e^{i\rho(s)}|\widehat f (s)|$, we call
$\rho(s)$ the \emph{phase} of $f$.  The phase is well defined
\emph{a.e.} modulo $2\pi$-translation.

\subsection{The Commutant of
$\{D,T\}$}
We have $\mathcal{F}\{D,T\}'\mathcal{F}^{-1} = \{\widehat
D,\widehat T\}'$. It turns out that $\{\widehat D,\widehat T\}'$ has
an elementary characterization in terms of Fourier multipliers:

\begin{theorem}\label{thm6}
\[
\{\widehat D,\widehat T\}' = \{M_h\colon \ h\in L^\infty(\mathbb{R}) \text{ and
} h(s) = h(2s) \text{ a.e.}\}.
\]
\end{theorem}

\begin{proof}
Since $\widehat D = D^*$ and $D$ is unitary, it is clear that $M_h \in
\{\widehat D, \widehat T\}'$ if and only if $M_h$ commutes with $D$. So let
$g\in L^2(\mathbb{R})$ be arbitrary. Then (a.e.) we have
\begin{align*}
(M_hDg)(s) &= h(s) (\sqrt 2\ g(2s)),\quad \text{and}\\
(DM_hg)(s) &= D(h(s) g(s)) = \sqrt h(2s) g(2s).
\end{align*}
Since these must be equal a.e.\ for arbitrary $g$, we must have $h(s) = h(2s)$
a.e.
\end{proof}

Now let $E = [-2,-1)\cup [1,2)$, and for $n\in \mathbb{Z}$ let $E_n=
\{2^nx\colon \ x\in E\}$. Observe that the sets $E_n$ are disjoint and have
union $\mathbb{R}\backslash\{0\}$. So if $g$ is any uniformly bounded function
on
$E$, then $g$ extends uniquely (a.e.) to a function $\tilde g\in
L^\infty(\mathbb{R})$ satisfying
\[
\tilde g(s) = \tilde g(2s),\qquad s\in \mathbb{R},
\]
by setting
\[
\tilde g(2^ns) = g(s),\qquad s\in E, n\in \mathbb{Z},
\]
and $\tilde g(0)=0$. We have $\|\tilde g\|_\infty = \|g\|_\infty$. Conversely,
if $h$ is any  function satisfying $h(s) = h(2s)$ a.e., then $h$ is uniquely
(a.e.) determined by its restriction to $E$. This 1-1 mapping $g\to M_{\tilde
g}$ from $L^\infty(E)$ onto $\{\widehat D, \widehat T\}'$ is a $*$-isomorphism.

We will refer to a function $h$ satisfying $h(s) = h(2s)$ a.e.\ as a
2-{\em dilation periodic function}. This gives a simple algorithm
for computing a large class of wavelets from a given one, by simply
modifying the \emph{phase} (see also section 4.7):
\begin{align}
&\text{\em Given $\psi$,, let $\widehat \psi = \mathcal{F}(\psi)$, choose a
real-valued function } h\in L^\infty(E)\nonumber\\
\label{eq22}
&\text{\em arbitrarily, let $g = \exp(ih)$, extend to a 2-dilation periodic}\\
&\text{\em function $\tilde g$ as above, and compute } \psi_{\tilde g} =
\mathcal{F}^{-1}(\tilde g\widehat\psi).\nonumber
\end{align}

In the description above, the set $E$ could clearly be replaced with
$[-2\pi,-\pi)\cup [\pi,2\pi)$, or with any other ``dyadic'' set $[-2a,a)\cup
[a,2a)$ for some $a>0$.

\subsection{The Shannon Wavelet}

We now give an account of $s$-elementary and $MSF$-wavelets. The two
most elementary dyadic orthonormal wavelets are the well-known {\em
Haar wavelet\/} and {\em Shannon's wavelet\/} (also called the
Littlewood--Paley wavelet).  The Shannon set is the prototype of the
class of wavelet sets.

Shannon's wavelet is the $L^2(\mathbb{R})$-function with Fourier transform
$\widehat\psi_S = \frac1{\sqrt{2\pi}} \chi_{E_0}$ where
\begin{equation}\label{eq33}
E_0 = [-2\pi, -\pi) \cup [\pi,2\pi).
\end{equation}
The argument that $\widehat\psi_S$ is a wavelet is in a way even more
transparent than for the Haar wavelet. And it has the advantage of generalizing
nicely. For a simple argument, start from the fact that the exponents
\[
\{e^{i\ell s}\colon \ n\in \mathbb{Z}\}
\]
restricted to $[0,2\pi]$ and normalized by $\frac1{\sqrt{2\pi}}$ is an
orthonormal basis for $L^2[0,2\pi]$. Write $E_0 = E_-\cup E_+$ where $E_- =
[-2\pi, -\pi)$, $E_+ = [\pi,2\pi)$. Since $\{E_- +2\pi, E_+\}$ is a partition of
$[0,2\pi)$ and since the exponentials $e^{i\ell s}$ are invariant under
translation by $2\pi$, it follows that
\begin{equation}\label{eq34}
\left\{\frac{e^{i\ell s}}{\sqrt{2\pi}}\Big|_{E_0}\colon \ n\in
\mathbb{Z}\right\}
\end{equation}
is an orthonormal basis for $L^2(E_0)$. Since $\widehat T = M_{e^{-is}}$, this
set can be written
\begin{equation}\label{eq35}
\{\widehat T^\ell \widehat\psi_s\colon \ \ell\in \mathbb{Z}\}.
\end{equation}
Next, note that any ``dyadic interval'' of the form $J = [b,2b)$, for some $b>0$
has the property that $\{2^nJ\colon \ n\in\mathbb{Z}\}$, is a partition of
$(0,\infty)$. Similarly, any set of the form
\begin{equation}\label{eq36}
\mathcal{K} = [-2a,-a)\cup [b,2b)
\end{equation}
for $a,b>0$, has the property that
\[
\{2^n\mathcal{K}\colon \ n\in \mathbb{Z}\}
\]
is a partition of $\mathbb{R}\backslash\{0\}$. It follows that the space
$L^2(\mathcal{K})$, considered as a subspace of $L^2(\mathbb{R})$, is a complete
wandering subspace for the dilation unitary $(Df)(s) = \sqrt 2\ f(2s)$. For each
$n\in \mathbb{Z}$,
\begin{equation}\label{eq37}
D^n(L^2(\mathcal{K})) = L^2(2^{-n}\mathcal{K}).
\end{equation}
So $\bigoplus_n D^n(L^2(\mathcal{K}))$ is a direct sum decomposition of
$L^2(\mathbb{R})$. In particular $E_0$ has this property. So
\begin{equation}\label{eq38}
D^n\left\{\frac{e^{i\ell s}}{\sqrt{2\pi}}\Big|_{E_0}\colon \ \ell\in
\mathbb{Z}\right\} = \left\{\frac{e^{2^ni\ell s}}{\sqrt{2\pi}}\Big|_{2^{-n}E_0}
\colon \ \ell\in \mathbb{Z}\right\}
\end{equation}
is an orthonormal basis for $L^2(2^{-n}E_0)$ for each $n$. It follows that
\[
\{D^n\widehat T^\ell \widehat\psi_s\colon \ n,\ell\in \mathbb{Z}\}
\]
is an orthonormal basis for $L^2(\mathbb{R})$. Hence $\{D^nT^\ell\psi_s\colon \
n,\ell\in \mathbb{Z}\}$ is an orthonormal basis for $L^2(\mathbb{R})$, as
required.

For our work, in order to proceed with developing an
operator-algebraic theory that had a chance of directly impacting
concrete function-theoretic wavelet theory we needed a large supply
of examples of wavelets which were elementary enough to work with.
First, we found another ``Shannon-type'' wavelet in the literature.
This was the Journe wavelet, which we found described on p.~136 in
Daubechies book [8]. Its Fourier transform is $\widehat \psi_J =
\frac1{\sqrt{2\pi}} \chi_{E_J}$, where
\[
E_J = \left[-\frac{32\pi}7, -4\pi\right) \cup \left[-\pi,
-\frac{4\pi}7\right)\cup \left[\frac{4\pi}7, \pi\right) \cup \left[4\pi,
\frac{32\pi}7\right).
\]
Then, thinking the old adage ``where there's smoke there's fire!'', we
painstakingly worked out many more examples. So far, these are the basic
building
blocks in the {\em concrete\/} part of our theory. By this we mean the part of
our theory that has had some type of direct impact on function-theoretic wavelet
theory.

\subsection{Definition of Wavelet Set} We define a {\em wavelet
set\/} to be a measurable subset $E$ of $\mathbb{R}$ for which
$\frac1{\sqrt{2\pi}} \chi_E$ is the Fourier transform of a wavelet.
The wavelet $\widehat\psi_E := \frac1{\sqrt{2\pi}}\chi_E$ is called
$s$-{\em elementary\/} in [5].

It turns out that this class of wavelets was also discovered and
systematically explored completely independently, and in about the
same time period, by Guido Weiss (Washington University), his
colleague and former student E. Hernandez (U.\ Madrid), and his
students X.\ Fang and X.\ Wang.  Two of the papers of this group are
[9] and [17], in which they are called MSF (minimally supported
frequency) wavelets.  In signal processing, the parameter $s$, which
is the independent variable for $\widehat\psi$, is the {\em
frequency\/} variable, and the variable $t$, which is the
independent variable for $\psi$, is the {\em time\/} variable. It is
not hard to show that no function with support a subset of a wavelet
set $E$ of strictly smaller measure can be the Fourier transform of
a wavelet. (Here, the support of a measurable function is defined to
be the set of points at which it does not vanish.) In other words,
an MSF wavelet has \emph{minimal} possible support in the frequency
domain.  However, the problem of whether the support set of any
wavelet necessarily contains a wavelet set remains open.    It was
raised by this author (Larson) in a talk about ten years ago, and
has been open for several years. We include it as Problem 3 in the
final section of this article.  A natural subproblem, which was
posed in the same talk, asks whether a wavelet with minimal possible
support in the frequency domain is in fact an MSF wavelet; or
equivalently, is its support a wavelet set?

\subsubsection{The Spectral Set Condition}



>From the argument above describing why Shannon's wavelet is, indeed, a wavelet,
it is clear that {\em sufficient\/} conditions for $E$ to be a wavelet set are

\begin{quote}
(i)~~the normalized exponential $\frac1{\sqrt{2\pi}} e^{i\ell s}$, $\ell\in
\mathbb{Z}$, when restricted to $E$ should constitute an orthonormal basis for
$L^2(E)$ (in other words $E$ is a {\em spectral set\/} for the integer lattice
$\mathbb{Z}$),
\end{quote}

\noindent and

\begin{quote}
(ii)~~The family $\{2^nE\colon \ n\in\mathbb{Z}\}$ of dilates of $E$ by integral
powers of 2 should constitute a measurable partition (i.e.\ a partition modulo
null sets) of $\mathbb{R}$.
\end{quote}

\noindent These conditions are also necessary. In fact if a set $E$ satisfies
(i), then for it to be a wavelet set it is obvious that (ii) must be satisfied.
To show that (i) must be satisfied by a wavelet set $E$, consider the vectors
\[
\widehat D^n \widehat\psi_E = \frac1{\sqrt{2\pi}} \chi_{2^{-n}E},\qquad n\in
\mathbb{Z}.
\]
Since $\widehat\psi_E$ is a wavelet these must be orthogonal, and so the sets
$\{2^nE\colon \ n\in~\mathbb{Z}\}$ must be disjoint modulo null sets. It follows
that $\{\frac1{\sqrt{2\pi}} e^{i\ell s}|_E\colon \ \ell\in \mathbb{Z}\}$ is not
only an orthonormal set of vectors in $L^2(E)$, it must also {\em span\/}
$L^2(E)$.

It is known from the theory of \emph{spectral sets} (as an
elementary special case) that a measurable set $E$ satisfies (i) if
and only if it is a generator of a measurable partition of
$\mathbb{R}$ under translation by $2\pi$ (i.e.\ iff $\{E+2\pi
n\colon \ n\in \mathbb{Z}\}$ is a measurable partition of
$\mathbb{R}$). This result generalizes to spectral sets for the
integral lattice in $\mathbb{R}^n$. For this elementary special case
a direct proof is not hard.

\subsection{Translation and Dilation Congruence}

We say that measurable sets $E,F$ are {\em translation congruent modulo\/}
$2\pi$ if there is a measurable bijection $\phi\colon \  E\to F$ such that
$\phi(s)-s$ is an integral multiple of $2\pi$ for each $s\in E$; or
equivalently, if there is a measurable partition $\{E_n\colon \ n\in
\mathbb{Z}\}$ of $E$ such that
\begin{equation}\label{eq39}
\{E_n  + 2n\pi\colon \ n\in \mathbb{Z}\}
\end{equation}
is a measurable partition of $F$. Analogously, define measurable sets $G$ and
$H$
to be {\em dilation congruent modulo\/} 2 if there is a measurable bijection
$\tau\colon \ G\to H$ such that for each $s\in G$ there is an integer $n$,
depending on $s$, such that $\tau(s) = 2^ns$; or equivalently, if there is a
measurable partition $\{G_n\}^\infty_{-\infty}$ of $G$ such that
\begin{equation}\label{eq40}
\{2^nG\}^\infty_{-\infty}
\end{equation}
is a measurable partition of $H$. (Translation and dilation congruency modulo
other positive numbers of course make sense as well.)

The following lemma is useful.

\begin{lemma}\label{lem7}
Let $f\in L^2(\mathbb{R})$, and let $E = \text{\rm supp}(f)$. Then $f$ has the
property that
\[
\{e^{ins}f\colon \ n\in \mathbb{Z}\}
\]
is an orthonormal basis for $L^2(E)$ if and only if
\begin{itemize}
\item[(i)] $E$ is congruent to $[0,2\pi)$ modulo $2\pi$, and
\item[(ii)] $|f(s)| = \frac1{\sqrt{2\pi}}$ a.e.\ on $E$.
\end{itemize}
\end{lemma}

If $E$ is a measurable set which is $2\pi$-translation congruent to $[0,2\pi)$,
then since
\[
\left\{\frac{e^{i\ell s}}{\sqrt{2\pi}}\Big|_{[0,2\pi)}\colon \ \ell\in
\mathbb{Z}\right\}
\]
is an orthonormal basis for $L^2[0,2\pi]$ and the exponentials $e^{i\ell s}$ are
$2\pi$-invariant, as in the case of Shannon's wavelet it follows that
\[
\left\{\frac{e^{i\ell s}}{\sqrt{2\pi}}\Big|_E\colon \ \ell\in \mathbb{Z}\right\}
\]
is an orthonormal basis for $L^2(E)$. Also, if $E$ is $2\pi$-translation
congruent to $[0,2\pi)$, then since
\[
\{[0,2\pi) + 2\pi n\colon \ n\in \mathbb{Z}\}
\]
is a measurable partition of $\mathbb{R}$, so is
\[
\{E + 2\pi n\colon \ n\in \mathbb{Z}\}.
\]
These arguments can be reversed.

We say that a measurable subset $G\subset \mathbb{R}$ is a 2-{\em dilation
generator\/} of a {\em partition\/} of $\mathbb{R}$ if the sets
\begin{equation}\label{eq41}
2^nG := \{2^ns\colon \ s\in G\},\qquad n\in \mathbb{Z}
\end{equation}
are disjoint and $\mathbb{R}\backslash \cup_n 2^nG$ is a null set. Also, we say
that $E\subset \mathbb{R}$ is a $2\pi$-{\em translation generator of a
partition\/} of $\mathbb{R}$ if the sets
\begin{equation}\label{eq42}
E + 2n\pi := \{s + 2 n\pi\colon \ s\in E\},\qquad n\in \mathbb{Z},
\end{equation}
are disjoint and $\mathbb{R}\backslash \cup_n (E+2n\pi)$ is a null set.

\begin{lemma}\label{lem8}
A measurable set $E\subseteq \mathbb{R}$ is a $2\pi$-translation generator of a
partition of $\mathbb{R}$ if and only if, modulo a null set, $E$ is translation
congruent to $[0,2\pi)$ modulo $2\pi$. Also, a measurable set $G\subseteq
\mathbb{R}$ is a 2-dilation generator of a partition of $\mathbb{R}$ if and only
if, modulo a null set, $G$ is a dilation congruent modulo 2 to the set $[-2\pi,
-\pi) \cup [\pi,2\pi)$.
\end{lemma}

\subsection{A Criterion} The following is a useful criterion for
wavelet sets. It was published independently by Dai--Larson in [5]
and by Fang--Wang (who were students of Guido Weiss) in [9] at about
the same time (in December, 1994). In fact, it is amusing that the
two papers had been submitted within two days of each other; only
much later did we even learn of each other's work on wavelets and of
this incredible timing.

\begin{proposition}\label{pro9}
Let $E\subseteq\mathbb{R}$ be a measurable set. Then $E$ is a wavelet set if and
only if $E$ is both a 2-dilation generator of a partition (modulo null sets) of
$\mathbb{R}$ and a $2\pi$-translation generator of a partition (modulo null
sets) of $\mathbb{R}$. Equivalently, $E$ is a wavelet set if and only if $E$ is
both translation congruent to $[0,2\pi)$ modulo $2\pi$ and dilation congruent to
$[-2\pi,-\pi) \cup [\pi,2\pi)$ modulo 2.
\end{proposition}

Note that a set is $2\pi$-translation congruent to $[0,2\pi)$ iff it is
$2\pi$-translation congruent to $[-2\pi, \pi)\cup [\pi,2\pi)$. So the last
sentence of Proposition \ref{pro9} can be stated:\ A measurable set $E$ is a
wavelet set if and only if it is both $2\pi$-translation and 2-dilation
congruent to the Littlewood--Paley set $[-2\pi, -\pi)\cup [\pi,2
\pi)$.

\subsection{Phases} If $E$ is a wavelet set, and if $f(s)$ is any
function with support $E$ which has constant modulus
$\frac1{\sqrt{2\pi}}$ on $E$, then $\mathcal{F}^{-1}(f)$ is a
wavelet. Indeed, by Lemma \ref{lem7} $\{\widehat T^\ell f\colon \
\in \mathbb{Z}\}$ is an orthonormal basis for $L^2(E)$, and since
the sets $2^nE$ partition $\mathbb{R}$, so $L^2(E)$ is a complete
wandering subspace for $\widehat D$, it follows that $\{\widehat
D^n\widehat T^\ell f\colon\ n,\ell\in \mathbb{Z}\}$ must be an
orthonormal basis for $L^2(\mathbb{R})$, as required. In [9, 17] the
term MSF-wavelet includes this type of wavelet. So MSF-wavelets can
have arbitrary phase and $s$-elementary wavelets have phase 0.  If
$\psi$ is a wavelet we say [5] that a real-valued function $\rho(s)$
is \emph{attainable} as a phase of $\psi$ if the function
$e^{i\rho(s)}|\psi(s)|$ is also the Fourier transform of a wavelet.
So \emph{every} phase is \emph{attainable} in this sense for an MSF
or $s$-elementary wavelet. Attainable phases of wavelets have been
studied in [5] and [26], in particular.

\subsection{Some Examples of One-Dimensional Wavelet Sets} It is
usually easy to determine, using the dilation-translation criteria,
in Proposition \ref{pro9}, whether a given finite union of intervals
is a wavelet set. In fact, to verify that a given ``candidate'' set
$E$ is a wavelet set, it is clear from the above discussion and
criteria that it suffices to do two things.

\begin{quote}
(1)~~Show, by appropriate partitioning, that $E$ is 2-dilation-congruent to a
set of the form $[-2a,-a)\cup [b,2b)$ for some $a,b>0$.
\end{quote}
and
\begin{quote}
(2)~~Show, by appropriate partitioning, that $E$ is $2\pi$-translation-congruent
to a set of the form $[c,c+2\pi)$ for some real number $c$.
\end{quote}

On the other hand, wavelet sets suitable for testing hypotheses can
be quite difficult to construct. There are very few ``recipes'' for
wavelet sets, as it were. Many families of such sets have been
constructed for reasons including perspective, experimentation,
testing hypotheses, etc., including perhaps the pure enjoyment of
doing the computations $-$ which are somewhat ``puzzle-like'' in
nature. In working with the theory it is nice (and in fact we find
it necessary) to have a large supply of wavelets on hand that permit
relatively simple analysis.

For this reason we take the opportunity here to present for the
reader a collection of such sets, mainly taken from [5], leaving
most of the work in verifying that they are indeed wavelet sets to
the reader.

We refer the reader to [6] for a proof of the existence of wavelet
sets in $\mathbb{R}^{(n)}$, and a proof that there are sufficiently
many to generate the Borel structure of $\mathbb{R}^{(n)}$. These
results are true for arbitrary expansive dilation factors. Some
concrete examples in the plane were subsequently obtained by Soardi
and Weiland, and others were obtained by Gu and by Speegle in their
thesis work at A\&M. Two had also been obtained by Dai for inclusion
in the revised concluding remarks section of our Memoir [5].

In these examples we will usually write intervals as half-open
intervals $[~\cdot ~,~)$ because it is easier to verify the
translation and dilation congruency relations (1) and (2) above when
wavelet sets are written thus, even though in actuality the
relations need only hold modulo null sets.

(i)~~As mentioned above, an example due to Journe of a wavelet which admits no
multiresolution analysis is the $s$-elementary wavelet with wavelet set
\[
\left[-\frac{32\pi}7, -4\pi\right)\cup \left[-\pi, \frac{4\pi}7\right)\cup
\left[\frac{4\pi}7, \pi\right) \cup \left[4\pi, \frac{32\pi}7\right).
\]
To see that this satisfies the criteria, label these intervals, in order, as
$J_1, J_2, J_3, J_4$ and write $J =\cup J_i$. Then
\[
J_1\cup 4J_2 \cup 4J_3\cup J_4 = \left[-\frac{32\pi}7, -\frac{16\pi}7\right)
\cup \left[\frac{16\pi}7, \frac{32\pi}7\right).\]
This has the form $[-2a,a)\cup [b,2b)$ so is a 2-dilation generator of a
partition of $\mathbb{R}\backslash\{0\}$. Then also observe that
\[
\{J_1 + 6\pi, J_2 +2\pi, J_3, J_4-4\pi\}
\]
is a partition of $[0,2\pi)$.

(ii)~~The Shannon (or Littlewood--Paley) set can be generalized. For
any $-\pi < \alpha < \pi$, the set
\[
E_\alpha = [-2\pi + 2\alpha, -\pi + \alpha) \cup [\pi + \alpha, 2\pi +2\alpha)
\]
is a wavelet set. Indeed, it is clearly a 2-dilation generator of a partition of
$\mathbb{R}\backslash\{0\}$, and to see that it satisfies the translation
congruency criterion for $-\pi < \alpha \le 0$ (the case $0<\alpha<\pi$ is
analogous) just observe that
\[
\{[-2\pi + 2\alpha, 2\pi) + 4\pi, [-2\pi, -\pi+\alpha) + 2\pi, [\pi +
\alpha,2\pi + 2\alpha)\}
\]
is a partition of $[0,2\pi)$. It is clear that $\psi_{E_\alpha}$ is then a
continuous (in $L^2(\mathbb{R})$-norm) path of $s$-elementary wavelets. Note
that
\[
\lim_{\alpha\to\pi} \widehat\psi_{E_\alpha} = \frac1{\sqrt{2\pi}} \chi_{[2\pi,
4\pi)}.
\]
This is {\em not\/} the Fourier transform of a wavelet because the set $[2\pi,
4\pi)$ is not a 2-dilation generator of a partition of
$\mathbb{R}\backslash\{0\}$. So
\[
\lim_{\alpha\to\pi} \psi_{E_\alpha}
\]
is not an orthogonal wavelet. (It is what is known as a Hardy wavelet because it
generates an orthonormal basis for $H^2(\mathbb{R})$ under dilation and
translation.) This example demonstrates that $\mathcal{W}(D,T)$ is {\em not\/}
closed in $L^2(\mathbb{R})$.

(iii)~~Journe's example above can be extended to a path. For $-\frac\pi7 \le
\beta\le \frac\pi7$ the set
\[
J_\beta = \left[-\frac{32\pi}7, -4\pi + 4\beta\right) \cup \left[-\pi +\beta,
-\frac{4\pi}7\right) \cup \left[\frac{4\pi}7, \pi+\beta\right)\cup \left[4 \pi +
4\beta, 4\pi + \frac{4\pi}7\right)
\]
is a wavelet set. The same argument in (i) establishes dilation congruency. For
translation, the argument in (i) shows congruency to $[4\beta, 2\pi+4\beta)$
which is in turn congruent to $[0,2\pi)$ as required. Observe that here, as
opposed to in (ii) above, the limit of $\psi_{J_\beta}$ as $\beta$ approaches
the boundary point $\frac\pi7$ {\em is\/} a wavelet. Its wavelet set is a union
of 3 disjoint intervals.

(iv)~~Let $A\subseteq [\pi, \frac{3\pi}2)$ be an arbitrary measurable subset.
Then there is a wavelet set $W$, such that $W\cap [\pi, \frac{3\pi}2)=A$. For
the
construction, let
\begin{align*}
B &= [2\pi, 3\pi)\backslash 2A,\\
C &= \left[-\pi, -\frac\pi2\right)\backslash (A-2\pi)\\
\text{and}\quad D &= 2A-4\pi.
\end{align*}
Let
\[
W = \left[\frac{3\pi}2, 2\pi\right)\cup A\cup B\cup C\cup D.
\]
We have $W\cap [\pi, \frac{3\pi}2) =A$. Observe that the sets
$[\frac{3\pi}2,2\pi)$, $A,B,C,D$, are disjoint. Also observe that the sets
\[
\left[\frac{3\pi}2,2\pi\right), A, \frac12 B, 2C, D,
\]
are disjoint and have union $[-2\pi,-\pi)\cup [\pi,2\pi)$. In addition, observe
that the sets
\[
\left[\frac{3\pi}2,2\pi\right), A,B-2\pi, C + 2\pi, D+2\pi,\]
are disjoint and have union $[0,2\pi)$. Hence $W$ is a wavelet set.

(v)~~Wavelet sets for arbitrary (not necessarily integral) dilation factors
other then 2 exist. For instance, if $d\ge 2$ is arbitrary, let
\begin{align*}
A &= \left[-\frac{2d\pi}{d+1}, -\frac{2\pi}{d+1}\right),\\
B &= \left[\frac{2\pi}{d^2-1}, \frac{2\pi}{d+1}\right),\\
C &= \left[\frac{2d\pi}{d+1}, \frac{2d^2\pi}{d^2-1}\right)
\end{align*}
and let $G = A \cup B\cup C$. Then $G$ is $d$-wavelet set. To see this, note
that
$\{A+2\pi,B,C\}$ is a partition of an interval of length $2\pi$. So $G$ is
$2\pi$-translation-congruent to $[0,2\pi)$. Also, $\{A,B,d^{-1}C\}$ is a
partition of the set $[-d\alpha, -\alpha)\cup [\beta,d\beta)$ for $\alpha =
\frac{2\pi}{d^2-1}$, and $\beta = \frac{2\pi}{d^2-1}$, so from this form it
follows that $\{d^nG\colon \ n\in \mathbb{Z}\}$ is a partition of
$\mathbb{R}\backslash\{0\}$. Hence if $\psi :=
\mathcal{F}^{-1}(\frac1{\sqrt{2\pi}} \chi_G)$, it follows that
$\{d^{\frac{n}2}\psi(d^nt-\ell)\colon \ n,\ell\in \mathbb{Z}\}$ is orthonormal
basis for $L^2(\mathbb{R})$, as required.

\section{Operator-Theoretic Interpolation for Wavelet Sets}

Operator-theoretic interpolation takes a particularly natural form
for the special case of s-elementary (or MSF) wavelets that
facilitates hands-on computational techniques in investigating its
properties. Let $E, F$ be a pair of wavelet sets. Then for (a.e.) $x
\in E$ there is a unique $y \in F$ such that $x - y \in
2\pi\mathbb{Z}$. This is the $\emph{translation congruence}$
property of wavelet sets.  Also, for (a.e.) $x \in E$ there is a
unique $z \in F$ such that $\frac{x}{z}$ is an integral power of
$2$. This is the \emph{dilation congruence} property of wavelet
sets. (See section 2.5.6.)

There is a natural \emph{closed-form algorithm} for the
$\emph{interpolation unitary}$ $V_{\psi_E}^{\psi_F}$ which maps the
wavelet basis for  $\widehat{\psi}_E $ to the wavelet basis for
$\widehat{\psi}_F$.  Indeed, using both the translation and dilation
congruence properties of $\{E,F\}$, one can explicitly compute a
(unique) measure-preserving transformation $\sigma := {\sigma}_E^F$
mapping $\mathbb{R}$ onto $\mathbb{R}$ which has the property that
$V_{\psi_E}^{\psi_F}$ is identical with the \emph{composition
operator} defined by:
$$ f \mapsto f \circ {\sigma}^{-1} $$
for all $f \in L^2(\mathbb{R})$. With this formulation, compositions
of the maps $\sigma$ between different pairs of wavelet sets are not
difficult to compute, and thus products of the corresponding
interpolation unitaries can be computed in terms of them.

\subsection{The Interpolation Map $\sigma$}
Let $E$ and $F$ be arbitrary wavelet sets. Let $\sigma\colon \ E\to
F$ be the 1-1, onto map implementing the $2\pi$-translation
congruence. Since $E$ and $F$ both generated partitions of
$\mathbb{R}\backslash\{0\}$ under dilation by powers of 2, we may
extend $\sigma$ to a 1-1 map of $\mathbb{R}$ onto $\mathbb{R}$ by
defining $\sigma(0)=0$, and
\begin{equation}\label{eq44}
\sigma(s) = 2^n\sigma(2^{-n}s) \quad \text{for}\quad s\in 2^nE,\quad
n\in \mathbb{Z}.
\end{equation}
We adopt the notation $\sigma^F_E$ for this, and call it the {\em
interpolation map\/} for the ordered pair $(E,F)$.

\begin{lemma}\label{lem5.1}
In the above notation, $\sigma^F_E$ is a measure-preserving
transformation from $\mathbb{R}$ onto $\mathbb{R}$.
\end{lemma}

\begin{proof}
Let $\sigma := \sigma^F_E$. Let $\Omega\subseteq \mathbb{R}$ be a
measurable set. Let $\Omega_n= \Omega \cap 2^nE$, $n\ni\mathbb{Z}$,
and let $E_n = 2^{-n} \Omega_n\subseteq E$. Then $\{\Omega_n\}$ is a
partition of $\Omega$, and we have $m(\sigma(E_n)) = m(E_n)$ because
the restriction of $\sigma$ to $E$ is measure-preserving. So
\begin{align*}
m(\sigma(\Omega)) &= \sum_n m(\sigma(\Omega_n)) = \sum_n m(2^n \sigma(E_n))\\
&= \sum_n 2^nm(\sigma(E_n)) = \sum_n 2^nm(E_n)\\
&= \sum_n m(2^nE_n) =\sum_n m(\Omega_n) = m(\Omega).
\end{align*}
\vskip-3em
\end{proof}

A function $f\colon \mathbb{R}\to \mathbb{R}$ is called 2-{\em
homogeneous\/} if $f(2s) = 2f(s)$ for all $s\in \mathbb{R}$.
Equivalently, $f$ is 2-homogeneous iff $f(2^ns) = 2^nf(s)$,
$s\in\mathbb{R}$, $n\in\mathbb{Z}$. Such a function is completely
determined by its values on any subset of $\mathbb{R}$ which
generates a partition of $\mathbb{R}\backslash\{0\}$ by 2-dilation.
So $\sigma^F_E$ is the (unique) 2-homogeneous extension of the
$2\pi$-transition congruence $E\to F$. The set of all 2-homogeneous
measure-preserving transformations of $\mathbb{R}$ clearly forms a
group under composition. Also, the composition of a
2-dilation-periodic function $f$ with a 2-homogeneous function $g$
is (in either order) 2-dilation periodic. We have $f(g(2s)) =
f(2g(s)) = f(g(s))$ and  $g(f(2s)) = g(f(s))$. These facts will be
useful.

\subsubsection{An Algorithm For The Interpolation Unitary}

 Now let
\begin{equation}\label{eq45}
U^F_E := U_{\sigma^F_E},
\end{equation}
where if $\sigma$ is any measure-preserving transformation of
$\mathbb{R}$ then $U_\sigma$ denotes the composition operator
defined by $U_\sigma f = f\circ\sigma^{-1}$, $f\in L^2(\mathbb{R})$.
Clearly $(\sigma^F_E)^{-1} = \sigma^E_F$ and $(U^F_E)^* = U^E_F$. We
have $U^F_E\widehat\psi_E = \widehat\psi_F$ since $\sigma^F_E(E)=F$.
That is,
\[
U^F_E\widehat\psi_E = \widehat\psi_E \circ\sigma^E_F =
\frac1{\sqrt{2\pi}} \chi_{_E} \circ \sigma^E_F = \frac1{\sqrt{2\pi}}
\chi_{_F} = \widehat\psi_F.
\]

\begin{proposition}\label{pro5.2}
Let $E$ and $F$ be arbitrary wavelet sets. Then $U^F_E\in
\mathcal{C}_{\widehat\psi_E}(\widehat D,\widehat T)$.  Hence
$\mathcal{F}^{-1} U^F_E \mathcal{F}$ is the interpolation unitary
for the ordered pair $(\psi_E,\psi_F)$.
\end{proposition}

\begin{proof}
Write $\sigma = \sigma^F_E$ and $U_\sigma = U^F_E$. We have
$U_\sigma\widehat\psi_E = \widehat\psi_F$ since $\sigma(E) = F$. We
must show
\[
U_\sigma\widehat D^n \widehat T^l \widehat\psi_E = \widehat
D^n\widehat T^l U_\sigma \widehat\psi_E,\quad n,l\in \mathbb{Z}.
\]
We have
\begin{align*}
(U_\sigma\widehat D^n\widehat T^l \widehat\psi_E)(s) &=
(U_\sigma\widehat D^n
e^{-ils} \widehat\psi_E)(s)\\
&= U_\sigma2^{-\frac{n}2} e^{-il2^{-n}s} \widehat\psi_E(2^{-n}s)\\
&= 2^{-\frac{n}2} e^{-il2^{-n}\sigma^{-1}(s)}
\widehat\psi_E(2^{-n}\sigma^{-1}(s))\\
&= 2^{-\frac{n}2} e^{-il\sigma^{-1}(2^{-n}s)}
\widehat\psi_E(\sigma^{-1}(2^{-n}s))\\
&= 2^{-\frac{n}2} e^{-il\sigma^{-1}(2^{-n}s)} \widehat\psi(2^{-n}s).
\end{align*}
This last term is nonzero iff $2^{-n}s\in F$, in which case
$\sigma^{-1}(2^{-n}s) = \sigma^E_F(2^{-n}s)$ $=2^{-n}s+2\pi k$ for
some $k\in \mathbb{Z}$ since $\sigma^E_F$ is a
$2\pi$-translation-congruence on $F$. It follows that
$e^{-il\sigma^{-2}(2^{-n}s)} = e^{-il2^{-n}s}$. Hence we have
\begin{align*}
(U_\sigma\widehat D^n\widehat T^l\widehat\psi_E)(s) &=
2^{-\frac{n}2}
e^{-ils^{-2n}s} \widehat \psi_F(2^{-n}s)\\
&= (\widehat D^n\widehat T^l\widehat\psi_F)(s)\\
&= (\widehat D^n\widehat T^lU_\sigma\widehat \psi_E)(s).
\end{align*}
We have shown $U^F_E \in \mathcal{C}_{\widehat\psi_E}(\widehat
D,\widehat T)$. Since $U^F_E\widehat\psi_E = \widehat\psi_F$, the
uniqueness part of Proposition ~3.1 shows that
$\mathcal{F}^{-1}U^F_E\mathcal{F}$ must be the interpolation unitary
for $(\psi_E,\psi_F)$.
\end{proof}

\subsection{The Interpolation Unitary Normalizes The Commutant}

\begin{proposition}\label{pro5.3}
Let $E$ and $F$ be arbitrary wavelet sets. Then the interpolation
unitary for the ordered pair $(\psi_E,\psi_F)$ normalizes
$\{D,T\}'$.
\end{proposition}

\begin{proof}
By Proposition \ref{pro5.2} we may work with $U^F_E$ in the Fourier
transform domain. By Theorem 6, the generic element of $\{\widehat
D,\widehat T\}'$ has the form $M_h$ for some 2-dilation-periodic
function $h\in L^\infty(\mathbb{R})$. Write $\sigma = \sigma^F_E$
and $U_\sigma= U^F_E$. Then
\begin{equation}\label{eq46}
U^{-1}_\sigma M_hU_\sigma= M_{h\circ\sigma^{-1}}.
\end{equation}
So since the composition of a 2-dilation-periodic function with a
2-homogeneous function is 2-dilation-periodic, the proof is
complete.
\end{proof}

\subsubsection{$\mathcal{C}_\psi(D,T)$ is Nonabelian}

It can also be shown ([5, Theorem 5.2 (iii)]) that if $E,F$ are
wavelet sets with $E\ne F$ then $U^F_E$ is not contained in the
double commutant $\{\widehat D,\widehat T\}''$. So since $U^F_E$ and
$\{\widehat D,\widehat T\}'$ are both contained in the local
commutant of $\mathcal{U}_{\widehat D,\widehat T}$ at
$\widehat\psi_E$, this proves that
$\mathcal{C}_{\widehat\psi_E}(\widehat D,\widehat T)$ is nonabelian.
In fact (see [5, Proposition 1.8]) this can be used to show that
$\mathcal{C}_\psi(D,T)$ is nonabelian for every wavelet $\psi$. We
suspected this, but we could not prove it until we discovered the
``right'' way of doing the needed computation using $s$-elementary
wavelets.

The above shows that a pair $(E,F)$ of wavelets sets (or, rather,
their corresponding $s$-elementary wavelets) admits
operator-theoretic interpolation if and only if Group$\{U^F_E\}$ is
contained in the local commutant
$\mathcal{C}_{\widehat\psi_E}(\widehat D,\widehat T)$, since the
requirement that $U^F_E$ normalizes $\{\widehat D,\widehat T\}'$ is
automatically satisfied. It is easy to see that this is equivalent
to the condition that for each $n\in\mathbb{Z}$, $\sigma^n$ is a
$2\pi$-congruence of $E$ in the sense that
$(\sigma^n(s)-s)/2\pi\in\mathbb{Z}$ for all $s\in E$, which in turn
implies that $\sigma^n(E)$ is a wavelet set for all $n$. Here
$\sigma = \sigma^F_E$. This property hold trivially if $\sigma$ is
{\em involutive\/} (i.e.\ $\sigma^2=$ identity).

\subsubsection{The Coefficient Criterion}
In cases where ``torsion'' is present, so $(\sigma^F_E)^k$ is the
identity map for some finite integer $k$, the von Neumann algebra
generated by $\{\widehat D,\widehat T\}'$ and $U :=U^F_E$ has the
simple form
\[
\left\{\sum^k_{n=0} M_{h_n}U^n\colon \ h_n\in L^\infty(\mathbb{R})
\text{ with } h_n(2s) = h_n(s),\quad s\in \mathbb{R}\right\},
\]
and so each member of this ``interpolated'' family of wavelets has
the form
\begin{equation}\label{eq47}
\frac1{\sqrt{2\pi}} \sum^k_{n=0} h_n(s) \chi_{\sigma^n(E)}
\end{equation}
for 2-dilation periodic ``coefficient'' functions $\{h_n(s)\}$ which
satisfy the necessary and sufficient condition that the operator
\begin{equation}\label{eq48}
\sum^k_{n=0} M_{h_n }U^n
\end{equation}
is unitary.

A standard computation shows that the map $\theta$ sending $\sum^k_0
M_{h_n}U^n$ to the $k\times k$ function matrix $(h_{ij})$ given by
\begin{equation}\label{eq49}
h_{ij} = h_{\alpha(i,j)}\circ \sigma^{-i+1}
\end{equation}
where $\alpha(i,j) = (i+1)$ modulo $k$, is a $*$-isomorphism. This
matricial algebra is the cross-product of $\{D,T\}'$ by the
$*$-automorphism $ad(U^F_E)$ corresponding to conjugation with
$U^F_E$. For instance,  if $k=3$ then $\theta$ maps
\[
M_{h_1} + M_{h_2} U^F_E + M_{h_3}(U^F_E)^2
\]
to
\begin{equation}\label{eq50}
\left(\begin{matrix}
h_1&h_2&h_3\\
h_3\circ\sigma^{-1}&h_1\circ\sigma^{-1}&h_2\circ\sigma^{-1}\\
h_2\circ\sigma^{-2}&h_3\circ\sigma^{-2}&h_1\circ\sigma^{-2}
\end{matrix}\right).
\end{equation}
This shows that $\sum^k_0 M_{h_n}U^n$ is a unitary operator iff the
scalar matrix $(h_{ij})(s)$ is unitary for almost all $s\in
\mathbb{R}$. Unitarity of this matrix-valued function is called the
{\em Coefficient Criterion\/} in [5], and the functions $h_i$ are
called the interpolation coefficients. This leads to formulas for
families of wavelets which are new to wavelet theory.

\subsection{Interpolation Pairs of Wavelet Sets}

 For many
interesting cases of note, the interpolation map  $\sigma^F_E$ will
in fact be an \emph{involution} of $\mathbb{R}$ (i.e. $\sigma \circ
\sigma = id$, where $\sigma := \sigma^F_E$, and where $id$ denotes
the identity map). So torsion \emph{will} be present, as in the
above section, and it will be present in an essentially simple form.
The corresponding interpolation unitary will be a \emph{symmetry} in
this case (i.e. a selfadjoint unitary operator with square $I$).

It is curious to note that verifying a simple operator equation $U^2
= I$ directly by matricial computation can be extremely difficult.
It is much more computationally feasible to verify an equation such
as this by pointwise (a.e.) verifying explicitly the relation
$\sigma \circ \sigma = id$ for the interpolation map.  In [5] we
gave a number of examples of interpolation pairs of wavelet sets. We
give below a collection of examples that has not been previously
published: Every pair sets from the Journe family is an
interpolation pair.

\subsection{Journe Family Interpolation Pairs}

Consider the parameterized path of \emph{generalized Journe} wavelet
sets given in Section ~4.8 Item (iii).  We have
$$ J_{\beta} = \left[-\frac{32\pi}{7} , -4\pi - 4\beta\right) \cup \left[-\pi +
\beta
, -\frac{4\pi}{7}\right) \cup\left[\frac{4\pi}{7}, \pi + \beta\right) \cup
\left[4\pi +
4\beta, 4\pi + \frac{4\pi}{7}\right)$$
where the set of parameters $\beta$ ranges $-\frac{\pi}{7} \leq \beta
\leq \frac{\pi}{7}$.

\begin{proposition} Every pair $(J_{\beta_1}, J_{\beta_2})$ is an
interpolation pair.
\end{proposition}
\begin{proof}
Let $\beta_1, \beta_2 \in \left[-\frac{\pi}{7}, \frac{\pi}{7}\right)$ with
$\beta_1
< \beta_2$. Write $\sigma = \sigma_{J_{\beta_2}}^{J_{\beta_1}} .$ We
need to show that
\begin{equation}{\sigma}^2(x)=x \tag{*}
\end{equation} for all $x \in \mathbb{R}$. Since $\sigma$ is
2-homogeneous, it suffices to verify (*) only for $x \in
J_{\beta_1}$. For $x \in J_{\beta_1} \cap J_{\beta_2}$ we have
$\sigma (x) = x$, hence ${\sigma}^2(x)=x$.  So we only need to check
(*) for $x \in (J_{\beta_1} \backslash J_{\beta_2})$. We have
$$J_{\beta_1} \backslash J_{\beta_2} = [-\pi + \beta_1 , -\pi +
\beta_2) \cup [4\pi + 4\beta_1 , 4\pi + 4\beta_2) .$$ It is useful
to also write
$$ J_{\beta_2} \backslash J_{\beta_1} = [-4\pi + 4\beta_1 , -4\pi +
4\beta_2) \cup [\pi + \beta_1 , \pi + \beta_2).$$

On $[-\pi + \beta_1 , -\pi + \beta_2 )$ we have $\sigma (x) = x +
2\pi$, which lies in $[\pi + \beta_1 , \pi + \beta_2)$.  If we
multiply this by $4$, we obtain
$4\sigma(x)\in[4\pi+4\beta_1,4\pi+4\beta_2) \subset J_{\beta_1}$.
And on $[4\pi + 4\beta_1, 4\pi + 4\beta_2)$ we clearly have $\sigma
(x) = x - 8\pi$, which lies in $[-4\pi + 4\beta_1 , -4\pi +
4\beta_2)$.

So for $x \in [-\pi + \beta_1 , -\pi + \beta_2)$ we have
$${\sigma}^2(x) = \sigma (\sigma (x) ) = \frac14 \sigma (4\sigma
(x)) = \frac14 [4\sigma(x) - 8\pi] = \sigma(x) -2\pi = x + 2\pi
- 2\pi = x.$$

On $[4\pi + 4\beta_1, 4\pi +4\beta_2)$ we have $\sigma(x) = x -
8\pi$, which lies in $[-4\pi + 4\beta_1 , -4\pi + 4\beta_2)$. So $\frac14 \sigma
(x) \in [-\pi + \beta_1, -\pi+\beta_2)$. Hence
$$\sigma\left(\frac14\sigma(x)\right) = \frac14\sigma(x) + 2\pi$$
and
thus
$${\sigma}^2(x) = 4\sigma\left(\frac14\sigma(x)\right) =
4\left[\frac14\sigma(x) + 2\pi\right] = \sigma(x) + 8\pi = x - 8\pi + 8\pi =
x$$ as required.

We have shown that for all  $x \in J_{\beta_1}$ we have
${\sigma}^2(x) = x$.  This proves that $(J_{\beta_1}, J_{\beta_2})$
is an interpolation pair.

\end{proof}

\section{Some Open Problems}

We will discuss four problems on wavelets that we have investigated
from an operator-theoretic point of view over the past ten years,
together with some related problems.   The set
 of orthonormal dyadic one-dimensional wavelets
is a set of vectors in the unit sphere of a Hilbert space $H =
L^2(\mathbb{R})$. It is natural to ask what are the topological
properties of $\mathcal {W}(D, T)$ as a subset of the metric space
$H$?  This type of question is interesting from a pure mathematical
point of view, speaking as an operator theorist, and it just may
have some practical consequences depending on the nature and the
degree of depth of solutions.

{\bf Problem 1:} [Connectedness]   This was the first global problem
in wavelets that we considered from an operator-theoretic point of
view. In [5] we posed a number of open problems in the context of
the memoir. The first problem we discussed was Problem A (in [5])
which conjectured that $\mathcal{W}(D, T)$ is
norm-arcwise-connected. It turned out that this conjecture was also
formulated independently by Guido Weiss and his group (see [9], [17,
[26]) from a harmonic analysis point of view (our point of view was
purely functional analysis), and this problem (and related problems)
was the primary stimulation for the creation of the WUTAM CONSORTIUM
-- a team of 14 researchers based at Washington University and Texas
A\&M University.  (See [26], for the first publication of this
group.) This \emph{connectedness conjecture} was answered {yes} in
[26] for the special case of the family of dyadic orthonormal MRA
wavelets in $L^2(\mathbb{R})$, but still remains open for the family
of \emph{arbitrary} dyadic orthonormal wavelets in
$L^2(\mathbb{R})$, as well as for the family of orthonormal wavelets
for any fixed $n$ and any dilation matrix in $\mathbb{R}^n$.  A
natural related problem which also remains open, is whether the set
of Riesz wavelets is connected.  An intermediate problem, which is
also open, asks whether given two orthonormal wavelets is there a
continuous path connecting them consisting of Riesz wavelets? (Some
evidence for a positive answer to this problem is given by
Proposition 6.1 below, which easily shows that every point on the
convex path connecting two wavelets (i.e. $(1-t)\psi + t\eta$)is a
Riesz wavelet \emph{except} for perhaps the midpoint corresponding
to $t = 0.5$. Thus this problem has an easy positive solution for
many pairs of orthonormal wavelets, but no way has been found to get
around the midpoint obstruction to show that all pairs are
connected, perhaps by some exotic type of path.) A subproblem is the
same problem but for the set of \emph{frame-wavelets}
$\mathcal{F}(D,T)$ (now widely called \emph{framelets}). Is the set
of all frame-wavelets connected?; or more specifically--is the set
of all Parseval frame wavelets connected?  The reader can easily
deduce some \emph{frame versions} of Proposition 6.1 using
elementary spectral theory of operators which provides some
quick-and-easy partial results on paths of frames. These are
tantalizing, but the main problems still remain open.

All of these connectivity problems have counterparts for other
unitary systems. For wavelet systems, they remain open (to our
knowledge) for all dimensions $n$ and all expansive matrix dilation
factors. And for other systems, in particular in [14], we showed
that for a fixed choice of modulation and translation parameters
(necessarily, of course, with product $\leq 1$) the set of
Weyl-Heisenberg (or Gabor) frames is connected in this sense, and
also it is norm-dense in $L^2(\mathbb{R})$ in the sense of Problem 2
below. Although the Weil-Heisenberg (aka Gabor) unitary systems are
of a simpler operator-theoretic structure than the wavelet systems,
and this permits the use of some techniques which do not work so
well in the wavelet theory, even so this perhaps is another reason
to think that general connectedness results are possible within the
wavelet theory.

{\bf Problem 2:} [Density]  Is the set $\mathcal{RW}(D,T)$ of all
Riesz wavelets \emph{dense} in the norm topology in the Hilbert
space $L^2(\mathbb{R})$?  This was posed as a conjecture by Larson
in a talk in August 1996 in a NATO conference held in Samos, Greece.
It was posed in the same spirit as the connectivity problem above,
in the sense that it asks about the topological nature of
$\mathcal{RW}(D,T)$ as a subset of the metric space
$L^2(\mathbb{R})$.  Like the connectivity problem it is a
\emph{global} type of problem.  A positive answer might be useful
for applications if it could be given a some type of
\emph{quantitative interpretation}. A subproblem of this, which was
discussed in several subsequent talks, is the same density problem
but for the set of frame-wavelets $\mathcal{F}(D,T)$.  (Of course a
positive answer for wavelets would imply it for framelets.) Like the
connectivity problem, this problem makes sense for the family of
Riesz wavelets for any fixed $n$ and any dilation matrix in
$\mathbb{R}^n$.

One of the reasons for thinking that this conjecture may be positive
is the following result,  which we think is the most
\emph{elementary} application of operator-theoretic interpolation.
It is abstracted from Chapter ~1 of [5], although the form in [5] is
a bit different.

\begin{proposition} Let $\mathcal{U}$ be a unitary system on a
Hilbert space $H$.  If $\psi_1$ and $\psi_2$ are in
$\mathcal{W(U)}$, then $$\psi_1 + \lambda \psi_2 \in
\mathcal{RW(U)}$$ for all complex scalars $\lambda$ with $|\lambda|
\neq 1$.  More generally, if $\psi_1$ and $\psi_2$ are in
$\mathcal{RW(U})$ then there are positive constants \space $b > a
 > 0$ such that $\psi_1 + \lambda \psi_2 \in
\mathcal{RW(U)}$ for all $\lambda \in \mathbb{C}$ with either
$|\lambda| < a$ or with $|\lambda| > b$.

\end{proposition}

\begin{proof}  If $\psi_1, \psi_2 \in \mathcal{W(U)}$, let $V$ be
the unique unitary in $\mathcal{C}_{\psi_2}(\mathcal{U})$ given by
Proposition ~3.1 such that $V\psi_2 = \psi_1$.  Then $$\psi_1 +
\lambda \psi_2 = (V + \lambda I)\psi_2 .$$  Since $V$ is unitary,
$(V + \lambda I)$ is an invertible element of
$\mathcal{C}_{\psi_2}(\mathcal{U})$ if $|\lambda | \neq 1$, so the
first conclusion follows from Proposition ~3.2.  Now assume $\psi_1,
\psi_2 \in \mathcal{RW(U)}$.  Let $A$ be the unique invertible
element of $\mathcal{C}_{\psi_2}(\mathcal{U})$ such that $A \psi_2 =
\psi_1$, and write $\psi_1 + \lambda \psi_2 = (A + \lambda
I)\psi_2$.  Since $A$ is bounded and invertible there are $b > a
> 0$ such that $$\sigma (A) \subseteq \{z \in \mathbb{C}: a < |z| < b
\}$$ where $\sigma (A)$ denotes the spectrum of $A$ , and the same
argument applies.

\end{proof}

The above proposition indicates that Riesz wavelets are plentiful.
As mentioned above, by writing $(1-t)\psi_2 + t\psi_1 = ((1-t)V +
tI)\psi_1$, and using the fact that the local commutant is a linear
space so contains $(1-t)V + tI$, it follows that a convex
combination of orthonormal wavelets is a Riesz wavelet except
possibly for the mid-point corresponding to $t = 0.5$~.  So, if
$\psi$ and $\eta$ are orthonormal wavelets, the line in the vector
space $L^2(\mathbb{R})$ containing the pair $\psi, \eta$ is in the
norm closure of the set of Riesz wavelets. Operator theoretic
interpolation shows that more general linear combinations of finite
families of Riesz wavelets are very often Riesz wavelets.  And if
one considers a finite family of wavelet sets $\mathcal{C}$, with
union $\mathcal{S}$, then the restricted sets of Riesz wavelets, or
orthonormal wavelets, or Parseval framelets, or framelets, which are
restricted in the sense that they have their frequency support
contained in $\mathcal{S}$, is always connected if the family of
wavelet sets is an \emph{interpolation family}, and these restricted
sets of wavelets are very often connected even if the interpolation
family criterion fails for $\mathcal{C}$.  (In fact, it is a
conjecture that these restricted sets of wavelets are \emph{always}
connected.) Moreover, these restricted sets are \emph{dense} in
$L^2(S)$ (considered as a subspace of $L^2(\mathbb{R})$) if the
family $\mathcal{C}$ is an interpolation family, and it is yet
another conjecture that they are \emph{always} dense in $L^2(S)$. It
is also known (see [5] for instance) that the linear span
$\mathcal{W}(D,T)$ is dense in $L^2(\mathbb{R})$. So these facts
together suggest that it is probably true that the set of Riesz
wavelets $\mathcal{RW}(D,T)$ is dense in $L^2(\mathbb{R})$. However,
the problem remains open. Also, as mentioned above under the
\emph{connectivity} problem, one can easily deduce some \emph{frame
versions} of Proposition 6.1 using elementary spectral theory of
operators which provides some quick-and-easy partial results on
density of Riesz wavelets and framelets. However, the general
problem for Riesz wavelets remains open. (For the density problem
for framelets, we mention that Marcin Bownick has recently obtained
a significant positive result!  It appears that he has solved the
problem positively for \emph{framelets}.  But apparently it remains
open for Riesz wavelets.)

As with the connectedness problem, all of these problems have
counterparts for other unitary systems. For wavelet systems, they
remain open (to our knowledge) for all dimensions $n$ and all
expansive matrix dilation factors (except for the recent interesting
framelet density result of Bownick mentioned above). And for other
systems, as mentioned in the context of Problem 1, in [14] we showed
that for a fixed choice of modulation and translation parameters
with product $\leq 1$ the set of Weyl-Heisenberg (or Gabor) frames
is dense in $L^2(\mathbb{R})$ in this sense.  This is another reason
to think that general density results might be possible within the
wavelet theory.

{\bf Problem 3:} [Frequency Support]  (See section 4.4.) [Must the
support of the Fourier transform of a wavelet contain a wavelet
set?] This conjecture was posed about 10 years ago, by Larson, and
the problem still remains open for the case of dimension $1$ and
dilation factor $2$. It makes sense for any finite dimension $n$ and
any matrix dilation, and it apparently remains unsolved in any case.
It has been studied by several researchers,  and Z. Rzeszotnik, in
particular, has made some progress on the problem. A related
problem, also posed by Larson, (see section 4.4) asks whether a
wavelet which has minimal support in the frequency domain is
necessarily an MSF wavelet. (In other words:  Is a minimal support
set in the frequency domain necessarily a wavelet set).

{\bf Problem 4:}  [Normalization] (See section 3.3.1.) If
$\{\psi,\eta\}$ is a pair of dyadic orthonormal wavelets, does the
interpolation unitary $V^\eta_\psi$ normalize $\{D,T\}'$? As
mentioned above, the answer is yes if $\psi$ and $\eta$ are
$s$-elementary wavelets. This problem makes sense for orthonormal
wavelets in higher dimensions for matrix dilation factors, and for
other scalar dilations in one dimension. We know of no
counterexample for any of these cases. However, the problem might
just lie in the fact that in most cases, other than for wavelet
sets, we have no reasonable techniques for doing the computations.
This problem was also discussed in context in Section 3.3.1, and it
 could be the most important problem remaining in the direction of
further development of the unitary system approach to wavelet
theory.

\end{document}